\input amstex 
\hsize=30truecc
\documentstyle{amsppt} 
\baselineskip=12pt 
\pagewidth{30pc}
\pageheight{52pc}
\loadbold
\TagsOnRight
\nologo
\nopagenumbers
\NoBlackBoxes  

\define\na{\Bbb N}

\define\sq{sequence}

\topmatter  
 
\title Spontaneous clustering in theoretical and some empirical stationary processes\endtitle
\rightheadtext {Spontaneous clustering in stationary processes} 
\leftheadtext {T. Downarowicz, Y. Lacroix and Didier Leandri} 
\author 
T. Downarowicz, Y. Lacroix and D. L\'eandri
\endauthor 
\address 
Institute of Mathematics and Computer Science, Wroclaw University of
Technology, Wy\-brze$\dot{\text z}$e Wys\-pia{\'n}\-skie\-go 27, 
50@-370 Wroc{\l}aw, 
Poland
\endaddress
\date August 26, 2008\enddate
\email 
downar\@pwr.wroc.pl
\endemail
\thanks
This paper was written during the first author's visit at ISITV in 2007.
\endthanks 
\address 
Institut des Sciences de l'Ing\'enieur de Toulon et du Var,
Laboratoire Syt\`emes Navals Complexes, 
Avenue G. Pompidou, B.P. 56, 83162 La Valette du Var Cedex,
France
\endaddress
\email 
yves.lacroix\@univ-tln.fr, joba\@club-internet.fr
\endemail
\subjclass 
\endsubjclass
\keywords stationary random process, return time, hitting time, attracting, limit law, cluster, the law of series
\endkeywords

\abstract 
In a stationary ergodic process, clustering is defined as the tendency of events to appear in series of increased
frequency separated by longer breaks. Such behavior, contradicting the theoretical ``unbiased behavior'' with exponential 
distribution of the gaps between appearances, is commonly observed in experimental processes and often difficult to explain. 
In the last section we relate one such empirical example of clustering, in the area of marine technology. In the theoretical part of the 
paper we prove, using ergodic theory and the notion of category, that clustering (even very strong) is in fact typical for 
``rare events'' defined as long cylinder sets in processes generated by a finite partition of an arbitrary (infinite aperiodic) 
ergodic measure preserving transformation. 
\endabstract 

\endtopmatter 
\document

\heading Introduction: recurrence and clustering\endheading 
Let us consider an aperiodic ergodic dynamical system $(X,\Cal B,\mu ,T)$, where $(X,\Cal B,\mu)$ is a standard probability space, and $T:X\to X$ is measurable, almost surely $1-1$, preserves $\mu$, for which it is an ergodic transformation [W]. 

Let us next consider a measurable partition $\Cal P=\{P_1,\ldots ,P_n\}$ of the space $X$. To this partition we associate the natural symbolic factor of the system, using the coding map 
$$c:x\in X\mapsto i\in\{1,\ldots ,n\}\Leftrightarrow x\in P_i$$
which generates the factor map [W] 
$$
\matrix
\pi:&X&\to&\{1,\ldots ,n\}^{\Bbb Z}\hfill\\
&x&\mapsto&(c(T^nx))_{n\in\Bbb Z} .\hfill\\
\endmatrix
$$
Then if $\nu=\pi\mu$ and $\sigma:\{1,\ldots ,n\}^{\Bbb Z}\to\{1,\ldots ,n\}^{\Bbb Z}$ denotes the left shift map, 
we obtain the factor $\pi : (X,\Cal B,\mu ,T)\to (\{1,\ldots ,n\}^{\Bbb Z},\Cal C,\nu ,\sigma)$ which satisfies $\pi\circ T=\sigma\circ \pi$ (the $\sigma$-algebra $\Cal C$ is the obvious one). 

This standard procedure, based on the selection of a partition on $X$, transforms the initial abstract dynamical system into a symbolic dynamical system. The coding procedure is somewhat natural in that it produces the accessible dynamics through the selection of a finite partition of the space of observables $X$. 

On the resulting symbolic space we inherit of standard subsets of the phase space $\{1,\ldots ,n\}^{\Bbb Z}$, which we call {\it cylinders, or blocks}. A typical such, of length $p$, is obtained by selecting a pattern $w\in\{1,\ldots ,n\}^p$, and defining 
$$
[w]=\{y\in\{1,\ldots ,n\}^{\Bbb Z}:(y_0,\ldots ,y_{p-1})=(w_0,\ldots ,w_{p-1})\}. 
$$
When the length $p$ of the cylinder $[w]$ increases, of course, the measure $\nu ([w])$ tends to $0$, which results from aperiodicity of the dynamics. 

Whence long cylinder sets are prototypical rare events (i.e. small measure sets) in the symbolic dynamics we accessed by cross-ruling our initial space $X$ with the partition $\Cal P$. The system is ergodic whence recurrence to any positive measure subset of $X$ must occur, a.s., and happens along a typical trajectory with frequency equal to the measure of the set (this is the ergodic theorem). 

Our paper concerns the study of recurrence to rare events generated by partitions as above.  If we set for given $B\in\Cal B$ with $\mu (B)>0$, 
$$
\tau_B(x)=\min\{k\ge 1:T^kx\in B\},
$$
then $\tau_B$ is $\mu$-a.s. well defined, integer valued. When $x\in B$, it is called the return time of $x$ to $B$, otherwise it is called the entry time of $x$ to $B$. Kac's theorem [K] states that $\sum_{k\ge 1}k\mu (\{\tau_B=k\}\cap B) =1$, whence the random variable $\mu (B)\tau_B$ defined on the probability space induced on $B$ has expected value equal to $1$. We therefore call $\mu (B)\tau_B$ the \lq\lq normalised return time\rq\rq . 

For rare events such as cylinder sets the distribution of return times or entry times to such events, though by the ergodic theorem has on the overall frequential distribution, may change quite a lot depending on the dynamics of the system. What is oftenly looked after is weak convergence of such distributions as the measure of the sets $B$ (resp. lengths of the cyinders) shrink to $0$ (resp. go to $\infty$). This is because if such convergence holds then the limiting distribution provides information about recurrence to rare events in the dynamical system, as it is an approximation of distributions converging to it weakly. 

Such a limit distribution is called an asymptotic for return times. We denote by $\tilde F_B$ the distribution function of the normalised return time to $B$, and by $F_B$ the distribution function of the entry time to $B$, i.e. of the random variable $\mu(B)\tau_B$. Following the same lines, asymptotics for entry times are analogously defined. 

Possible asymptotics have been characterized only recently ([L] for return times, [K-L] for entry times). 
An integral formula connecting asymptotics for entry times and the one for return times has been provided in [H-L-V], where it is proved that weak convergence of distributions, whenever it holds, must hold simultaneously for entry and return times. If $\tilde F$ is the weak limit distribution function for return times and $F$ is the one for entry times, then additionally 
$$
F(t)=\int_0^t(1-\tilde F(s))ds,\; t\ge 0. \tag$\star$
$$

Many research papers have been devoted to the study of return times of specific events, so-called cylinder sets, in stationary ergodic processes. We refer also the reader to expository papers [C] and [A-G] for further information. 
Most asymptotics along cylinder sets were found in mixing enough dynamical systems [W], but in all cases were proved to be exponential with parameter one, the distribution of which is the only fixed point of $(\star)$. Furthermore, in the treated cases, entropy of the system was positive [W]. 

An essential progress in understanding some phenomena towards asymptotics in processes with positive entropy has been obtained recently in [D-L], where it is proved that whatever the system, as soon as its entropy is positive, then any asymptotic $F$ for entry times along cylinder sets must satisfy 
$$
F(t)\le 1-e^{-t}, \; t\ge 0.\tag$\star\star$ 
$$
This was interpreted as a first explanation to a complicated and misunderstood common-sense phenomenon known as \lq\lq{\it the law of series}\rq\rq : indeed in any ergodic system, if $t>0$ and $\mu (B)>0$, if we define the variable 
$$
I(x)=\#\{0\le n\le \frac t{\mu (B)}:T^nx\in B\}, 
$$
then invariance of the measure implies that 
$$
\Bbb E(I)\approx t,\tag$\star\star\star$
$$
with uniform accuracy [D]. Considering an independent symbolic process so as to be one for which most randomness occurs, and using the well-known fact that in such system the asymptotics along cylinder sets exist and are exponential, we are led to argue that a small cylinder set $B$ has a tendency to occur presenting clusters more frequently than it would in an independent process, if given that $B$ occurs, the conditional expectation of $I$ increases, compared to what it is in the independent case. 

That is to say, $B$ clusters (or appears \lq\lq in series\rq\rq ) in distribution, if 
$$
\Bbb E(I\vert I>0)\ge \Bbb E_{Ind}(I_{Ind}\vert I_{Ind}>0), 
$$
where the subscript $Ind$ refers to the independent process. Now using $(\star\star\star)$ and the fact that 
$\Bbb E(I\vert I>0)=\Bbb E(I)\bigm/\mu (I>0)$, observing moreover that $\mu (I>0)=F_B(t)$, and finally approximating 
$F_{B,Ind}(t)\approx 1-e^{-t}$, we deduce 
$$
B\text{ clusters}\Leftrightarrow F_B(t)\le 1-e^{-t}. 
$$

In [D-L] we introduced the notions of attracting, strong attracting, repelling, and neutral recurrences for a set $B$, in reference to the above explained comparison to the case of the independent process. Below is a purely naive illustration of the phenomena we are describing, interpreted as \lq\lq along a typical orbit\rq\rq . 
$$
\align 
\text{\eightpoint unbiased \ }& 
..{...|}........{||....|..|....|}.....{...|}......{|..|}..........|...||......
...|......|.|...{...|}...{...||}....|...|.\\
\text{\eightpoint attracting \ }&..{...|}.............\underline{||..|.|..|}........
{...|}.......\underline{||}.........\underline{|.||}..........|..........\underline{|..|.|...||..|}.......|..\\
\vspace {5pt}
\text{\eightpoint strong attr. \ }
&......\underline{|.||||.|}..............................................\underline{||.|||..|.|}..
.....................\underline{||.||.|.||}....\\
\endalign
$$ 

This is the interpretation along which it was stated in [D-L] that $(\star\star)$ can be understood differently, saying that in positive entropy systems, clustering for rare events must be at least what it reveals to be in the neutral independent case. In other words, \lq\lq laws of series\rq\rq\ (more frequent clustering for rare events) correspond to the natural behaviour for positive entropy systems. 

This was completed in [D-L] by the study of asymptotics along cylinder sets for varying partitions $\Cal P$ of the space $X$, so as to understand what generically (in a Baire category setting) could happen to be the case (the set of partitions can be turned to a Polish structured space). 

It was proved that generically (on $\Cal P$) there exists an upper density one set $\Cal N\subset \Bbb N$ of lengths of cylinder sets along which asymptotics exist, and for entry times, converge to the degenerated distribution function $F\equiv 0$, the one corresponding to strong attracting above (enormous clusters). 

In the present paper we go forward in this last direction and prove that this genericity holds without the positive entropy assumption. We also addressed the question to specialists, since the \lq\lq law of series\rq\rq , frequently referred to under more pessimistic interpretations like \lq\lq Murphy laws\rq\rq , is used by engineers and manufacturers and formalised under the label \lq\lq Murphy proof devices\rq\rq, in aeronautics for instance. In fact other occurrences of Murphy proof procedures appear, and we reveal such in the last part of this paper, were we relate to some rare but clustering disturbances that occur on low speed flights of underwater gliders, without further understanding at this point of knowledge. 

The only reasonable understanding for this generic clustering phenomenon we can talk out, at this point, more like a guess, is the observation that rare events, to occur, need favourable conditions, probably even rarer, and that therefore when they are collected, conditionally, the event has a better measure and comes through easier. 
 
\medskip
The paper is organised as follows : the next section introduces rigourous formulations about genericity. The next section presents the formulation and  proof of our main result (Theorem 1).
In the formulation we use the notion of strong
clustering explained above in the Introduction. The last section is a presentation of an empirical example of a strongly clustered process occurring in marine technology, preceded by a short description of the performed experiment and its technical background.

\heading Typicality of strong attraction without entropy assumptions\endheading

Let us go back to the original dynamical system $(X,\Cal B,\mu,T)$ on which the observed symbolic dynamics is defined by the partition $\Cal P$. 

Denote by $\goth P_l$ the collection of all measurable partitions $\Cal P$ of $X$ into $l$ cells. 
In general there is no canonical measure on the space $\goth P_l$. The meaning of ``almost surely'' with respect 
to a random partition does not have a definite meaning. Instead, we will adopt the topological approach, and 
``typicality'' defined in terms of so-called {\it category}.

Recall that in a complete metric space a subset is called {\it residual} if it contains a dense $G_\delta$ set, 
equivalently, if its complement if of first category (i.e., is a countable union of nowhere dense sets). The Baire 
category theorem asserts that the intersection of any countable collection of residual sets is still residual. 
Because similar property is enjoyed by sets of measure 1 in probability spaces, residual sets are considered a 
topological analog of sets of measure 1. Given a fixed residual set, its elements are referred to as ``typical''.

If $(X,\Cal B,\mu)$ is a probability space, the {\it Rokhlin metric} endows $\goth P_l$ with a structure 
of a complete metric space. The distance in this metric between two $l$-element partitions $\Cal P, \Cal Q$ 
is defined as
$$
d(\Cal P, \Cal Q) = \inf_\pi \sum_{A\in\Cal P}\mu(A\,\triangle\,\pi(A)),
$$
where $\pi$ ranges over all bijections from $\Cal P$ to $\Cal Q$ and $\triangle$ denotes the symmetric difference 
of sets. A partition $\Cal Q$ which is very close to $\Cal P$ in this distance may be considered a slight perturbation
of $\Cal P$. Our theorem below applies to processes typical in the following sense : given an ergodic dynamical system $(X,\Cal B,\mu,T)$ and some $l\in\Bbb N$, the set of all $l$-element partitions $\Cal P$ of $X$
such that the generated process $(\{1,\ldots ,\ell\}^{\Bbb Z},\Cal C,\nu_{\Cal P} ,\sigma)$ satisfies the assertion of 
the theorem is residual in the Rokhlin metric in $\goth P_l$.

We will also use the phrase that a property $\Phi$ holds for all blocks of ``majority'' of lengths. 
By this we mean that there exists a set $\Bbb N_\Phi$ of {\it upper density} 1 such that if $n\in \Bbb N_\Phi$
then all blocks of of length $n$ satisfy $\Phi$. Upper density is defined as 
$$
\overline D(\Bbb N_\Phi)=\limsup_{n\to\infty}\frac {\#(\Bbb N_\Phi\cap[1,n])}n.
$$

\heading Formulation of the Theorem and the proof \endheading
Below we give two formulations of the main result. The first one is short thanks to the terminology introduced above 
and appeals to the intuitive understanding of the subject. The latter formulation is the rigorous version without
using the shortcut terminology.

\proclaim{Theorem 1}
In a typical ergodic process, for majority of lengths $n$, all rare elementary events of length $n$ reveal strong clustering.
\endproclaim

\proclaim{Theorem 1}
Let $(X,\Cal B,\mu,T)$ be an ergodic not periodic dynamical system. Fix some natural $l\ge 2$. Then in the space $\goth P_l$ 
of all \,$l$-element measurable partitions of $X$ endowed with the Rokhlin metric there exists a residual subset $\goth C$ such 
that for every $\Cal P\in\goth C$ the generated process $(\{1,\ldots ,\ell\}^{\Bbb Z},\Cal C,\nu_{\Cal P} ,\sigma)$ has the 
following property $\Phi$: There exists a set $\Bbb N_0\subset \Bbb N$ of upper density one, such that for every $\epsilon>0$ 
there is $n_\epsilon\in\Bbb N$ such that for every $n\in\Bbb N_0$, $n>n_\epsilon$ and every block 
$B$ of length $n$, the $\tilde F_B(\varepsilon)<\epsilon^2$.
\endproclaim

\demo{Proof} 
Fix $\epsilon>0$ and $N\in\na$. Suppose a partition $\Cal P\in\goth P_l$ satisfies the following property $\Phi_{\epsilon,N}$ : 
\medskip 
\noindent\lq\lq For every $n\in[N,N^2]$ and every block $B$ of length $n$ holds $\tilde F_B(\varepsilon)<\epsilon^2$.\rq\rq  
\medskip 
\noindent Clearly, if we perturb 
$\Cal P$ very little, the property will still be satisfied. Thus $\Phi_{\epsilon,N}$ holds on an open set 
$\goth C_{\epsilon, N}$ of partitions. Of course, the set 
$$
\goth C_{\epsilon} = \bigcup_{N\ge 1}\goth C_{\epsilon, N},
$$ 
of partitions such that the same property holds for some $N$, is also open. The main effort in the proof will be to show 
that this set is also dense. Once this is done, the proof is complete, because then the dense $G_\delta$ set $\goth C$ of 
partitions which fulfill the hypothesis of the theorem is obtained by intersecting the sets $\goth C_{\epsilon}$ over 
countably many parameters $\epsilon$ converging to zero. Every element of this intersection satisfies $\Phi_{\epsilon,N}$ 
for arbitrarily small $\epsilon$ and some $N$ depending on $\epsilon$. Clearly, $N$ must grow to infinity as $\epsilon$
decreases to zero. Notice that for any infinite \sq\ of natural numbers $N$ the set $\bigcup[N,N^2]$ has upper density 
1 in $\na$. 

It remains to prove the density of the open set $\goth C_{\epsilon}$. The proof (as most proofs of typicality) may look a 
bit artificial, as it is done by perturbing an arbitrarily chosen partition $\Cal P$ in a very specific way, so that a 
highly particular partition is created. This perturbation is then shown to belong to the set $\goth C_{\epsilon}$. It 
is important to realize that once the density is proved, the ``largeness'' properties of open dense sets imply that 
$\goth C_{\epsilon}$ is represented in the vicinity of $\Cal P$ by many more partitions, not only the artificially 
constructed perturbation. 

We begin with a technical lemma.

\proclaim{Lemma 1}
In every ergodic not periodic dynamical system $(X,\Cal B,\mu,T)$, for each sufficiently large $r\in\na$ there exists a 
``semiperiodic $r$-marker'', i.e., a measurable set $F_r$ such that the return time to 
$F_r$ assumes almost surely only two values: $r$ and $r+1$. 
\endproclaim
\demo{Proof} 
We need a measurable and shift-invariant procedure dividing the trajectory of almost every point $x\in X$ into intervals
of the two lengths $r$ and $r+1$. (Shift-invariant means that the division of the trajectory of $Tx$ coincides with the
division of trajectory of $x$ with shifted enumeration.) Once this is done, the set $F_r$ is defined as the collection 
of all points which have a division marker at the coordinate zero. 

Using ergodicity and nonperiodicity (or the Rokhlin theorem, see [R]) it is very easy to construct a set of positive measure such 
that the return time to this set assumes almost surely only values larger than or equal to $r^2$. The times of visits to 
this set divide each trajectory into intervals of lengths at least $r^2$ (see also [L]). 
Every such interval can be further divided into intervals of two lengths $r$ and $r+1$, and we can fix one such 
way for every length $m\ge r^2$ (for example, we can choose the division which minimizes the number of longer intervals 
and all longer intervals appear to the right of the shorter ones). This (clearly measurable) technique divides almost 
every trajectory into desired pieces in a shift-invariant way, as required.
\qed\enddemo

We continue with the main proof; we are proving that $\goth C_\epsilon$ is dense. Fix any $l$-element partition $\Cal P$ 
and some $\delta>0$. We will construct a perturbation
$\Cal P'$ within the distance $\delta$ from $\Cal P$, which belongs to $\goth C_{\epsilon,N}$ for some $N$. 
Since the partition has at least two elements, we select two of them and label them 0 and 1. 
Pick $L$ so large such that there exist $K=\frac 2{\epsilon^2}$ different blocks $W'_k$ ($k=1,2,\dots,K$) of length 
$\frac{L-1}2$, none of them equal to $000\dots0$, each of measure at most $\frac{\delta}{2LK}$ (here we admit also 
blocks that do not appear in the system and hence have measure zero). Denote by $W_k$ the block 
$1W_k'1000\dots0$ (with $\frac{L-1}2$ zeros) of length $L$. The blocks $W_1,\dots,W_K$ also have measures at most 
$\frac{\delta}{2LK}$. Notice that any two (different or equal) blocks from this family never occur with an overlap.
Let $r$ be an integer larger than $\frac {2L}{\delta}$. Let $N$ be larger than $2r+2$. Every cylinder set of positive 
measure over a block 
$B$ of a length $n\ge N$ decomposes into a finite number of sets depending on the positioning of the $r$-markers 
(there is at least one such marker in every occurrence of $B$). We denote these sets by $B_i$ ($i=1,2,\dots,j(B)$). 
Let $M$ be so large that all points $x$, except in a set $Z$ whose measure is smaller than $\frac 13$, satisfy the 
following: For every block $B$ of positive measure and length between $N$ and $N^2$, and for every $i =1,2,\dots,j(B)$ 
the orbit of $x$ visits $B_i$ between times 0 and $M$ at least $\frac 3\epsilon$ times. Now let $r_1 = KM$.

At this point we modify the partition $\Cal P$ as follows: In the $\Cal P$-name of (almost) every $x$ the interval 
between two consecutive $r_1$-markers splits into $K$ intervals of length $M$ (the last one may be of length $M+1$). We call them 
``sectors''. In the sector number $k$ ($k=1,2,\dots,K$) we put the block $W_k$ immediately to the right of every $r$-marker 
(replacing whatever was there), and next we replace all occurrences of all blocks $W_j$ ($j=1,2,\dots,K$) also by $W_k$.  
Because the blocks $W_j$ do not overlap, such exchange happens on disjoint intervals. Notice that we have changed the 
partition $\Cal P$ only on the $r$-marker set $F_r$ and its $L$ consecutive preimages by $T$ (jointly a set of measure $\frac\delta2$), and on the cylinders corresponding to the blocks $W_1,\dots,W_K$, and also their $L$ preimages by $T$ 
(jointly another set of measure $\frac\delta2$). Thus the Rokhlin distance between $\Cal P$ and the modified partition 
$\Cal P'$ is at most $\delta$. 

Let $B'$ be an arbitrary block of some length $n$ between $N$ and $N^2$ appearing with positive probability in the process 
generated by the modified partition. Consider an $x\in X$. By ergodicity, $B'$ occurs almost surely in the $\Cal P'$-name 
of $x$ at some position $n$. The block $B'$ is long enough 
so that at least one $r$-marker occurs within its length. Next to this marker $B'$ contains one of the blocks $W_k$. 
Because $W_k$ does not appear in any sector other then the sector number $k$ (of an interval between two $r_1$-markers) 
$B'$ also can occur only in sectors that carry the number $k$. In the $\Cal P$-name of $x$ (i.e., before modification) 
at the position $n$ there is some block $B$ and also there is a fixed positioning of $r$-markers along this block,
i.e., $T^n(x)$ belongs to some $B_i$. We know that the orbit of $x$ visits $B_i$ at least $\frac 3\epsilon$ times 
in each sector, except the cases when at the beginning of the sector the orbit falls into $Z$ (which happen with probability $<\frac13$). 
If this sector happens to be the sector number $k$, in the modified partition every such visit generates an occurrence 
of $B'$. 
Thus we conclude the following: The block $B'$ occurs in the $\Cal P'$-name of every $x\in Z$ at least $\frac 3\epsilon$ 
times in more than $\frac23$ of all sectors number $k$ and zero times in sectors carrying other numbers 
(see Figure 2; stars indicate the $r_1$-markers). 

$$
...|......\overset{*}\to{|}\underset{s.\,1}\to{\vphantom{_p}.......}
|\underset{s.\,2}\to{\vphantom{_p}.......}|.......
|\underset{s.\,k}\to{\!_{^{B'B'\!B'}}\!}|.......|.......|.......|.......
|\underset{s.\,K}\to{\vphantom{_p}.......}
\overset{*}\to{|}\underset{s.\,1}\to{\vphantom{_p}.......}
|\underset{s.\,2}\to{\vphantom{_p}.......}|.......
|\underset{s.\,k}\to{\!_{^{B'B'\!B'}}\!}
|.......|...
$$
{\eightpoint\it Figure 2: The distribution of the occurrences of $B'$.}
\bigskip\noindent

It remains to compute the intensity, normalize the waiting time, and prove that $\tilde F_{B'}(\epsilon)\le \epsilon^2$.
The intensity of the occurrences of $B'$ is at least $\frac 23\frac 3{\epsilon KM} = \frac 2{\epsilon KM}$. 
(the multiplier $\frac 23$ takes into account the visits in $Z$). The waiting time for the signal is smaller 
than $M$ only when the zero coordinate falls in a sector number $k$ or in the preceding sector. This happens with 
probability at most $\frac 2K=\epsilon^2$. Otherwise the normalized waiting time is larger than 
$M\frac {2}{\epsilon KM}=\epsilon$. This ends the proof.
\qed\enddemo

\heading
Clustering in an empirical process
\endheading

One of the multiple illustrations of the law of series appears in the precise study of the behavior of an ``underwater glider'',
one of the newest inventions in marine technology.
The underwater glider is an autonomous robot able to cross thousands of miles without any propeller [B-D-H-L]. 
By changing the buoyancy of the ballasts, the vehicle tends to ``fly'' alternatively downwards or upwards to the surface. 
As the vehicle is equipped with wings - hence has some gliding ability, the variations of buoyancy are changed into 
horizontal motion. Very little energy is needed to modify the buoyancy and, as a consequence, the linear flight between 
the top and bottom points of the trajectory of the glider is costless in terms of energy (see Figure 3).

\bigskip
\input epsf.tex
\epsfxsize=8truecm
\centerline{\epsfbox{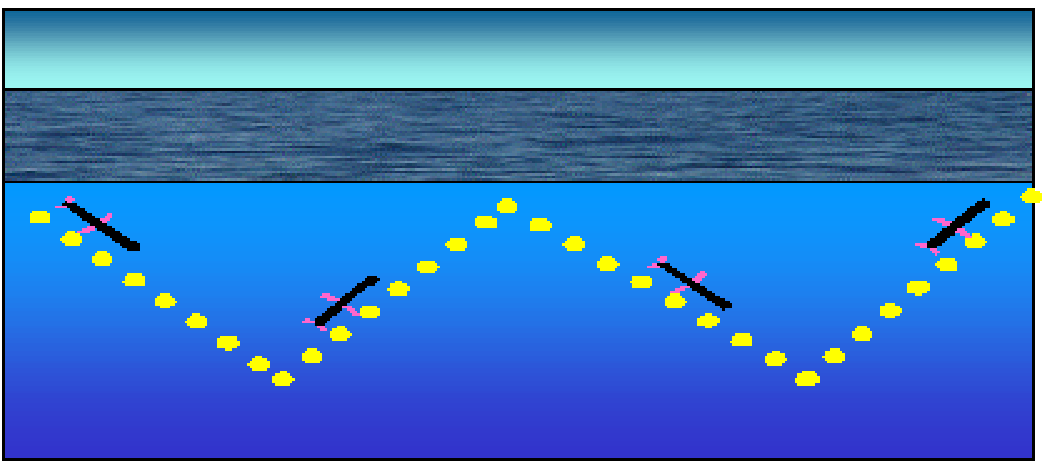}}
\smallskip
{\eightpoint\it Figure 3: The trajectory of a low energy underwater glider.}
\bigskip\noindent

This is the positive aspect of this brand new technology. However, during our numerous and lengthy tests at the sea 
we discovered some strange side effects, mainly linked to the fact that the glider  flies at a 
very low speed and therefore is very sensitive to hydrodynamic phenomena.

It was discovered [B-D-H-L] that even though one makes no doubt that the flow around the glider is linear, and the physical system (glider plus flow) is not chaotic (zero entropy), some heavy turbulences appear during short 
periods of time. They are encountered under a conjunction of circumstances with on ocean scale are pretty rare. The point is that it has been observed that these rare circumstances produce heavy turbulences, and those turbulences have a tendency to repeat overexpectedly. After a while, the stable flow returns. The authors decided to call theses turbulences ``germs'', due to some 
similarities with the fact that the cavitations originating the perturbations could be activated by some ``germs'' in the water, and \lq\lq pulsed out\rq\rq\ somewhat once at a time.

The authors were unable to predict, and still aren't, when and why these turbulences appear in a given homogeneous layer of water in the ocean. What they know by experience, is that once an unexpected turbulence appears, then others follow repeatedly with short gaps in time until the series stops for a new long period of smooth linear flow. This corresponds exactly to 
the pattern described as ``strong clustering''.

This series type behavior of appearances of turbulences has been integrated so as to produce, when the first turbulence of a series appears, some warning, and consequently for autonomous vehicles, the control law of the underwater vehicle has been reinforced in these flight specific phases.

The next step would be to produce some statistical inference, and derive from such risk estimation for stability of trajectory under reinforced control laws of different levels of security. This has of course immediate industrial consequences. But it is a difficult mathematical challenge.

\medskip 

\noindent http://www.labo-snc.eu/

\noindent http://www.pwr.wroc.pl/
\medskip 

\enddocument